\documentclass[12pt]{article}

\usepackage{bbm}
\usepackage{amssymb}
\usepackage{amsmath}
\usepackage{amsthm}
\usepackage{amsmath}
\usepackage{amsthm}
\usepackage{amsfonts}

\DeclareMathOperator{\1}{\mathbbm{1}}

\newcommand{\eee}{{\rm e}}

\newcommand{\mmp}{\mathbb{P}}

\newcommand{\dod}{\overset{{\rm d}}{\to}}
\newcommand{\tp}{\overset{\mmp}{\to}}

\newcommand{\me}{\mathbb{E}}

\theoremstyle{definition}

\theoremstyle{remark}

\begin{document}
\title{A comment on the article `The harmonic descent chain' by D.~J. Aldous, S. Janson and X. Li}
\date{\today}

\author{Alexander Iksanov\footnote{Faculty of Computer Science and Cybernetics, Taras Shevchenko National University of Kyiv, Ukraine, e-mail:
iksan@univ.kiev.ua}}

\maketitle

\begin{abstract}
An alternative proof is given for the main result of the article referred to in the title and published in ECP (2024). The proof exploits the theory of regenerative composition structures due to Gnedin and Pitman. The present article is a slight revision of my note written up in May 2024 as a response to the preprint arxiv.org version of the paper by Aldous, Janson and Li.  
\end{abstract}


\bigskip\bigskip

In \cite{Aldous et al:2024}, the authors investigate the discrete-time Markov chain $(X_n)_{n\geq 0}$ on positive integers with transition probabilities
$$p(j,j-i):=\frac{1}{ih_{j-1}},\quad 1\leq i < j,~ j \geq 2$$
and $p(1,1) = 1$. Here, $h_n:=\sum_{k=1}^n k^{-1}$ for $n\geq 1$. They state in Theorem 1.1 of \cite{Aldous et al:2024} and prove that, for each $i\geq 1$,
\begin{equation}\label{eq:main}
\lim_{n\to\infty}\mmp\{X_k=i+1~~\text{for some}~k|X_0=n+1\}=\frac{h_i}{\zeta(2)i},
\end{equation}
where $\zeta(2)=\sum_{k\geq 1}k^{-2}=\pi^2/6$.

The purpose of this note is to unveil and exploit a connection with an infinite `balls-in-boxes' scheme and also the theory of regenerative composition structures due to Gnedin and Pitman \cite{Gnedin+Pitman:2005}. 
One benefit of this finding is an effective proof of formula \eqref{eq:main}.

Let $S:=(S(t))_{t\geq 0}$ be a subordinator (an increasing L{\'e}vy process) with $S(0)=0$, zero drift, no killing and the L{\'e}vy measure $\nu$ defined by
\begin{equation*} 
\nu({\rm d}x)=\frac{\eee^{-x}}{1-\eee^{-x}}\1_{(0,\infty)}(x)\,{\rm d}x.
\end{equation*}
Let $E_1,\ldots, E_n$ be independent random variables with the exponential distribution of unit mean, which are independent of $S$. The closed range of $S$ has zero Lebesgue measure and splits the positive halfline into infinitely many disjoint intervals that we call gaps. We call a gap occupied if it contains at least one $E_j$, $j=1,2,\ldots, n$. To obtain the aforementioned `balls-in-boxes' scheme, take the gaps in the role of boxes and the points of the exponential sample in the role of balls. Put $X_0(n)=n$ and, for $k\geq 1$, denote by $X_k(n)$ the number of exponential points (out of $n$) which do not fall into the first $k$ occupied gaps counted left-to-right. According to Theorem 5.2(i) in \cite{Gnedin+Pitman:2005} which is an important result of the theory of regenerative composition structures developed by Gnedin and Pitman, $(X_k(n))_{k\geq 0}$ is a decreasing Markov chain which goes from $j$ to $j-i$ with probability
\begin{multline*}
\frac{{j \choose i}\int_{(0,\infty)}(1-\eee^{-x})^i \eee^{-x(j-i)}\nu({\rm d}x)}{\sum_{k=1}^j {j \choose i}\int_{(0,\infty)}(1-\eee^{-x})^i \eee^{-x(j-i)}\nu({\rm d}x)}=\frac{{j \choose i}\int_0^1 x^{i-1}(1-x)^{j-i}{\rm d}x}{\sum_{k=1}^j {j \choose i}\int_0^1 x^{i-1}(1-x)^{j-i}{\rm d}x}\\=\frac{1}{ih_j},\quad 1\leq i\leq j,\quad j\geq 1.
\end{multline*}
This formula entails the basic observation for the present note: for each $i\geq 1$ and $n>i$, $$\mmp\{X_k=i+1~~\text{for some}~k|X_0=n+1\}=\mmp\{X_k(n)=i~~\text{for some}~k\}.$$ Thus, \eqref{eq:main} follows once we have shown that
\begin{equation}\label{eq:main2}
\lim_{n\to\infty}\mmp\{X_k(n)=i~~\text{for some}~k\}=\frac{h_i}{\zeta(2)i},\quad i=1,2,\ldots
\end{equation}

\noindent {\sc Proof of \eqref{eq:main2}}. As usual, $\dod$ and $\tp$ will denote convergence in distribution and convergence in probability, respectively.

Put $S^\leftarrow(t):=\inf\{u\geq 0: S(u)>t\}$ for $t\geq 0$ and let $E_{1,n}<E_{2,n}<\ldots<E_{n,n}$ denote the order statistics of the sample $E_1,\ldots, E_n$. Then \begin{multline*}
\mmp\{X_k(n)=i~~\text{for some}~k\}=\mmp\{S(S^\leftarrow(E_{n-i,n}))<E_{n-i+1,n}\}\\=\mmp\{S(S^\leftarrow(E_{n-i,n}))-E_{n-i,n}<E_{n-i+1,n}-E_{n-i, n}\}.
\end{multline*}
The random variable $E_{n-i+1,\,n}-E_{n-i,\,n}$ has the exponential distribution of mean $1/i$ and is independent of $E_{n-i,\,n}$. Also,
$E_{n-i+1,\,n}-E_{n-i,\,n}$ is independent of $S$, hence of $S^\leftarrow$, because this is the case for the $E_k$ by assumption. Thus, $$\mmp\{X_k(n)=i~~\text{for some}~k\}=\me \big[\exp(-i (S(S^\leftarrow(E_{n-i,\,n}))-E_{n-i,\,n}))\big].$$ The subordinator $S$ is nonarithmetic with\footnote{More generally, according to a Hurwitz identity (see, for instance, formula (23.2.7) in \cite{Abram}), $$\int_0^\infty x^r \nu({\rm d}x)=\int_0^1 x^{-1}(-\log (1-x))^r{\rm d}x=\Gamma(r+1)\sum_{k\geq 1} k^{-r-1},\quad r>0,$$ where $\Gamma$ is the Euler gamma-function.} $\me [S(1)]=\int_{(0,\infty)}x\nu({\rm d}x)=\zeta(2)$. By Theorem 1 in \cite{Bertoin+Harn+Steutel:1999}, $S(S^\leftarrow(t))-t \dod \chi$ as $t\to\infty$, where $\chi$ is a random
variable with distribution $\mmp\{\chi>y\}=(\zeta(2))^{-1}\int_y^\infty \nu((x,\infty)){\rm d}x$ for $y\geq 0$. Since $E_{n-i,\,n}\tp +\infty$ as $n\to\infty$ and $E_{n-i,\,n}$ is independent of $(S(S^\leftarrow(t)))_{t\geq 0}$, we infer $S(S^\leftarrow(E_{n-i,\,n}))-E_{n-i,\,n}\dod \chi$ as $n\to\infty$. By the continuity theorem for Laplace transforms, the latter entails
\begin{multline*}
\lim_{n\to\infty} \me \big[\exp(-i (S(S^\leftarrow(E_{n-i,\,n}))-E_{n-i,\,n}))\big]=\me[\exp(-i\chi)]=\frac{1}{\zeta(2)}\int_0^\infty \eee^{-iy}\nu((y,\infty)){\rm d}y\\=\frac{1}{\zeta(2)i}\int_{(0,\infty)}(1-\eee^{-iy})\nu({\rm d}y)=
\frac{1}{\zeta(2)i}\int_0^\infty (1-\eee^{-iy})\frac{\eee^{-y}}{1-\eee^{-y}}{\rm d}y\\=\frac{1}{\zeta(2)i}\int_0^\infty \sum_{k=1}^i \eee^{-ky}{\rm d}y=\frac{h_i}{\zeta(2)i}.
\end{multline*}
The proof of \eqref{eq:main2} is complete.


\end{document}